\documentclass[11pt,amssymb]{amsart}
\usepackage{amssymb}
\usepackage[mathscr]{eucal}
\usepackage[all,cmtip]{xy}

\newcommand{\Z}{{\mathbb Z}}

\newcommand{\C}{{\mathbb C}}

\newcommand{\Q}{{\mathbb Q}}

\newcommand{\R}{{\mathbb R}}

\newtheorem{thm}{Theorem}[section]

\newtheorem{prop}[thm]{Proposition}

\newcommand{\Ks}{K^{\rm sep}}

.240pk scaled 1200 .240pk
\DeclareFontFamily{U}{wncy}{}
    \DeclareFontShape{U}{wncy}{m}{n}{<->wncyr10}{}
    \DeclareSymbolFont{mcy}{U}{wncy}{m}{n}
    \DeclareMathSymbol{\Sha}{\mathord}{mcy}{"58}
\begin{document}

\title[Recent developments]{Recent developments in the theory of linear algebraic groups: Good reduction and finiteness properties
}

\author[A.~Rapinchuk]{Andrei S. Rapinchuk}
\author[I.~Rapinchuk]{Igor A. Rapinchuk}

\address{Department of Mathematics, University of Virginia,
Charlottesville, VA 22904-4137, USA}

\email{asr3x@virginia.edu}

\address{Department of Mathematics, Michigan State University, East Lansing, MI
48824, USA}

\email{rapinchu@msu.edu}

\maketitle

\section*{}
The focus of this article is on the notion of good reduction for algebraic varieties and particularly for algebraic groups.
Among the most influential results in this direction is the Finiteness Theorem of G.~Faltings, which confirmed Shafarevich's conjecture that, over a given number field, there exist only finitely many isomorphism classes of abelian varieties of a given dimension that have good reduction at all valuations of the field lying outside a fixed finite set.
Until recently, however, similar questions have not been systematically considered in the context of linear algebraic groups.

Our goal  is to discuss the notion of good reduction for reductive linear algebraic groups and then formulate a finiteness conjecture concerning forms with good reduction over arbitrary finitely generated fields.
This conjecture has already been confirmed in a number of cases, including for algebraic tori over all fields of characteristic zero and also for some semi-simple groups, but the general case remains open. What makes this conjecture even more interesting are its multiple connections with other finiteness properties of algebraic groups, first and foremost, with the conjectural properness of the global-to-local map in Galois cohomology that arises in the investigation of the Hasse principle. Moreover, it turns out that techniques based on the consideration of good reduction have applications far beyond the theory of
algebraic groups: as an example, we will discuss a finiteness result arising in the analysis of length-commensurable Riemann surfaces that relies heavily on this approach. So, we hope that mathematicians working in different areas will find familiarity with good reduction quite rewarding and potentially useful. We therefore invite the reader to explore with us the fascinating symbiosis of ideas from the theory of algebraic groups and arithmetic geometry that lead to conjectures that are likely to be at the center of current efforts to develop the arithmetic theory of algebraic groups over fields more general than global.


\section{Reduction techniques in arithmetic geometry: a brief survey}

It has been known since antiquity that reduction techniques can be used to show the absence of  integral or rational solutions.
For example, taking the equation
$
x^2 - 7y^2 = -1
$
modulo 7, one easily finds that it does not have any integral solutions.
At the beginning of the 20th century, it became apparent that by reducing a given equation modulo various primes, one can not only detect the (non)-existence
of solutions, but in fact obtain information about the structure of the solution set. One of the earliest, and perhaps most revealing, examples arose in the study of elliptic curves. Let us consider an elliptic curve $E$ over the field of rational numbers $\mathbb{Q}$ given by its affine equation
\begin{equation}\label{E:1}
y^2 = f(x),
\end{equation}
where $f(x) = x^3 + ax + b$ is a cubic polynomial with integer coefficients without multiple roots, hence having {\it nonzero} discriminant $\Delta(f) = -4a^3 - 27b^2$. Recall that $E$ is the projective curve obtained by adding one point $O$ at infinity to the solution set of (\ref{E:1}), and that the chord-tangent law for addition of points makes the set $E(\mathbb{Q})$ of rational points into an abelian group.

In 1922, Louis Mordell showed that the group $E(\mathbb{Q})$ is finitely generated. One of the key steps in the argument is the so-called {\it Weak Mordell Theorem}, stating that the quotient $E(\mathbb{Q})/2E(\mathbb{Q})$ is finite.
The proof involves reducing the original equation (\ref{E:1}) modulo rational primes and considering those primes at which our curve has bad reduction. More precisely, reducing (\ref{E:1}) modulo a prime $p > 3$, we obtain the equation $y^2 = \bar{f}(x)$, where $\bar{f} = x^3 + \bar{a} x + \bar{b}$ and $\bar{a}$ and $\bar{b}$ denote the residues of $a$ and $b$ modulo $p$, respectively. We say that (1) has {\it good reduction} at $p$ if $\bar{f}$ does not have multiple roots, and {\it bad reduction} otherwise. Geometrically, at primes of good reduction, the reduced equation still defines an elliptic curve (over $\mathbb{F}_p = \mathbb{Z}/p\mathbb{Z}$), while at primes of bad reduction it defines a singular rational curve. (Informally, the reduced curve at primes of good reduction retains the ``type" of the original curve, while at primes of bad reduction it changes the type.) A somewhat subtle point here is that a change of coordinates does not essentially change an elliptic curve, i.e. results in an isomorphic curve, but may change the defining equation (\ref{E:1}). So, we say that an elliptic curve $E$ has {\it good reduction} at $p > 3$ if after a possible $\Q$-defined change of coordinates it can be given by an equation (\ref{E:1}) that has good reduction at $p$. (For example, the equation $y^2 = x^3 - 625 x$ has bad reduction at $p = 5$, but the elliptic curve it defines is isomorphic to the elliptic curve given by $y^2 = x^3 - x$, which has good reduction at $p = 5$.) Otherwise, $E$ has {\it bad reduction} at $p$.

For a given $E$, the set $S$ of primes of bad reduction is contained in the set of primes divisors of the discriminant $\Delta(f)$ for the equation (\ref{E:1}) defining $E$ in {\it any} coordinate system, and in particular is finite. In the proof  of the Weak Mordell Theorem, $S$ plays a crucial role. More specifically, in the case where $f$ has three rational roots, the proof yields an estimate of the order of $E(\mathbb{Q})/2E(\mathbb{Q})$, and hence of the rank of $E$, that depends only on the size of $S$. We refer the interested reader to \cite{Silverman} for the details.

One may  further wonder if the primes of bad reduction have an effect not only on the group of rational points, but on the elliptic curve $E$ itself. In his 1962 ICM talk, Shafarevich pointed out that if $S$ is a finite set of rational primes, then there are only {\it finitely many} isomorphism classes of elliptic curves over $\mathbb{Q}$ having good reduction at all $p \notin S$.


We note that the proof of the Weak Mordell Theorem can be extended to any number field and to higher-dimensional analogues of elliptic curves,  known as {\it abelian varieties}.
On the other hand, while the proof of Shafarevich's theorem can be extended to any number field $K$ and any finite set of places $S$,  it is very specific to elliptic curves. However, Shafarevich felt that his theorem was an instance of a far more general phenomenon, which prompted him to formulate the following finiteness conjecture for abelian varieties:

\vskip1mm

\noindent {\it Let $K$ be a number field and $S$ be a finite set of primes of $K$. Then for every $g \geq 1$, there exist only finitely many $K$-isomorphism classes of abelian varieties of dimension $g$ having good reduction at all primes $\mathfrak{p} \notin S$.}

This conjecture was proved by Faltings in 1982 as a culmination of research on finiteness properties in diophantine geometry over the course of several decades. Due to its numerous implications, connections, and generalizations (cf. \cite{Darmon}),
the subject remains one of the major themes in arithmetic geometry.

\section{An overview of algebraic groups}

To prepare for the discussion of good reduction of (linear) algebraic groups and their forms, we will now review some standard terminology.  The reader who is familiar with the theory of algebraic groups  may easily skip most of this section.

\subsection{Basic notions and examples} An {\it algebraic group} is subgroup $G \subset \mathrm{GL}_n(\Omega)$ of the general linear group over a ``large" algebraically closed field $\Omega$ (e.g., $\C$ in characteristic 0) defined by a set of polynomial conditions; in other words, $G$ is a subgroup that is closed in the Zariski topology on $\mathrm{GL}_n(\Omega)$.
Given algebraic groups $G \subset \mathrm{GL}_n(\Omega)$ and $H \subset \mathrm{GL}_m(\Omega)$, a {\it morphism} of algebraic groups $\varphi \colon G \to H$ is a group homomorphism that admits a ``polynomial" representation, i.e. there exist
\begin{equation}\label{E:AG1}
f_{kl} \in R_n(\Omega) := \Omega[x_{11}, \ldots , x_{nn}, (\det(x_{ij}))^{-1}]
\end{equation}
for $k,l = 1, \ldots , m$,
such that $\varphi(a) = (f_{kl}(a))$ for all $a \in G$. A morphism that admits an inverse morphism is called an {\it isomorphism}.

Next, let $K \subset \Omega$ be a subfield. 
When $K$ is {\it perfect}, an algebraic group $G \subset \mathrm{GL}_n(\Omega)$ is $K$-{\it defined} or is a $K$-{\it group} if it can be described by polynomial conditions with coefficients in $K$. Over general fields, this notion becomes a bit more technical: one needs to require that the ideal $I(G) \subset R_n(\Omega)$ of all functions that vanish on $G$ be generated by polynomials with coefficients in $K$.

A morphism $\varphi \colon G \to H$ between two $K$-groups is said to be $K$-{\it defined} if the $f_{kl}$ in (\ref{E:AG1}) can be chosen in $R_n(K) = K[x_{11}, \ldots , x_{nn}, (\det(x_{ij}))^{-1}]$. It may happen that for two $K$-groups $G$ and $H$, there exists an isomorphism $\varphi \colon G \to H$ defined over a field extension $F$ of $K$, but possibly not over $K$.
These issues are addressed in the theory of $F/K$-forms (see \S\ref{S-Forms}) and are of central importance since the structure of algebraic groups over algebraically closed fields is well understood (see \S2.2 for a brief outline and \cite{BorelAG} for a detailed account).

While this basic approach to (linear) algebraic groups is certainly viable, a key drawback is that
a given algebraic group $G \subset \mathrm{GL}_n(\Omega)$ may have a multitude of  embeddings $G \hookrightarrow \mathrm{GL}_m(\Omega)$, resulting in many different sets of defining equations. This obviously creates certain inconveniences in dealing with notions that may depend on the choice of these equations
(such as good reduction). The object that remains invariant for {\it all} realizations
is the $\Omega$-algebra of regular functions $\Omega[G] := R_n(\Omega)/I(G)$ in our previous notations. To keep track of the $K$-structure when $G$ is defined over a subfield $K \subset \Omega$, one needs to consider the algebra of $K$-regular functions $K[G] = R_n(K)/(I(G) \cap R_n(K))$. Without going into details, we note that the product operation on $G$ makes these algebras into commutative Hopf algebras, leading us to view $G$ as a group scheme over $\Omega$ and $K$, respectively. Although the perspective of the theory of group schemes is very natural for treating such topics as good reduction, in the present article we will minimize the use of this language.
So, for the most part we will work with concrete defining relations associated with a given matrix realization of $G$. One element that we will borrow from the schematic approach, though, is the possibility of avoiding the reference to a ``large" algebraically closed field $\Omega$, allowing us to write $G \subset \mathrm{GL}_n$ instead of $G \subset \mathrm{GL}_n(\Omega)$.

\vskip2mm

\noindent {\bf Example 2.1.} Standard examples include the {\it special linear group} $\mathrm{SL}_n = \{\, a \in \mathrm{GL}_n \ \vert \ \det(a) = 1 \}$ and  the {\it symplectic group} $\mathrm{Sp}_{2n} = \{\, a \in \mathrm{GL}_n \ \vert \ {}^ta E a = E \, \}$, where $E = \left( \begin{array}{cc} O_n & I_n \\
-I_n & O_n \end{array} \right)$ with $O_n$ and $I_n$ being, respectively, the zero and the identity $(n \times n)$-matrices. Both groups are defined over the prime subfield (i.e. $\Q$ in characteristic 0 and $\mathbb{F}_p$ in characteristic $p$).

Furthermore, let $q$ be a nondegenerate quadratic form with coefficients in a field $K$ of characteristic $\neq 2$, and let $Q$ be its matrix. Then the {\it orthogonal group} $\mathrm{O}_n(q) = \{\, a \in \mathrm{GL}_n \ \vert \ {}^ta Q a = Q \, \}$ and the {\it special orthogonal group} $\mathrm{SO}_n(q) = \mathrm{O}_n(q) \cap \mathrm{SL}_n$ are $K$-defined algebraic groups.

We will also encounter the {\it diagonal group} $\mathbb{D}_n$ and the {\it upper unitriangular group} $\mathbb{U}_n$ (both are again defined over the prime subfield).

\subsection{Structure of algebraic groups}
For any algebraic group $G$, the connected component of the identity $G^{\circ}$ for the Zariski topology is a normal subgroup of finite index (e.g., $\mathrm{SO}_n(q)$ is the connected component of $\mathrm{O}_n(q)$).  So, the analysis of arbitrary algebraic group essentially reduces to the analysis of connected ones. We will now introduce three disjoint classes of algebraic groups (unipotent groups, algebraic tori, and simple groups), and then describe how an arbitrary connected group can be constructed from these.

\vskip1mm

\noindent {\it Unipotent groups.} An algebraic group $U \subset \mathrm{GL}_n$ is called {\it unipotent} if the relation $(u - I_n)^n = O_n$ holds on it identically. Then there exists $g \in \mathrm{GL}_n(\Omega)$ such that $gUg^{-1} \subset \mathbb{U}_n$.

\vskip1mm

\noindent {\it Algebraic tori.} An algebraic group $T \subset \mathrm{GL}_n$ is {\it diagonalizable} if there exists $g \in \mathrm{GL}_n(\Omega)$ such that $gTg^{-1} \subset \mathbb{D}_n$. A connected diagonalizable group $T$ is called an {\it algebraic  torus}; it is isomorphic to $\mathbb{D}_m$, where $m = \dim T$.  If $T$ is defined over $K$ and an element $g$ conjugating $T$ into $\mathbb{D}_n$ can be chosen in $\mathrm{GL}_n(K)$, then $T$ is said to be $K$-{\it split}, and there is an isomorphism $T \simeq \mathbb{D}_m$ defined over $K$. A $K$-torus $T$ that does not contain any $K$-split subtori of positive dimension is called $K$-{\it anisotropic}.

\vskip1mm

\noindent {\it Simple and semi-simple groups.} A connected algebraic group $G$ is called {\it absolutely almost simple} if it is noncommutative and does not have any proper connected normal algebraic subgroups. Then $G$ automatically has a finite center $Z(G)$. To describe the classification of absolutely almost simple groups, we recall that a {\it central isogeny} $\pi \colon G \to H$ is a surjective morphism of connected algebraic groups whose schematic kernel is finite and central (in characteristic 0, this simply means that
the usual kernel is finite). An absolutely almost simple group $G$ is {\it simply connected} (resp., {\it adjoint}) if there is no nontrivial isogeny $\pi \colon H \to G$ (resp., $\pi \colon G \to H$).

Given an absolutely almost simple group $G$, one fixes a maximal torus $T \subset G$ and then considers the corresponding  root system $\Phi(G , T)$ (cf. \cite[Ch. IV, \S14]{BorelAG}) . The latter  turns out to be reduced and irreducible, hence belongs to one of the types $\textsf{A}_n, \textsf{B}_n, \ldots , \textsf{G}_2$. The main result is that over an an algebraically (and even separably) closed field, absolutely almost simple groups are classified up to isogeny by the type of the root system, and each isogeny class contains a unique (up to isomorphism) simply connected group $\widetilde{G}$ and a unique adjoint group $\overline{G}$. While the classification of general absolutely almost simple groups over a non-closed field $K$ is a very hard open problem (cf. \S5), the classification of $K$-{\it split groups} (i.e. those that contain a maximal $K$-torus that splits over $K$) is essentially the same as over the closed fields.

A group $G$ is {\it semi-simple} if there exists a central isogeny
$
\pi \colon G_1 \times \cdots \times G_r \to G,
$
where the $G_i$'s are absolutely almost simple group.

\vskip2mm

\noindent {\bf Example 2.2.} The groups $\mathrm{SL}_n$ and $\mathrm{Sp}_{2n}$ are absolutely almost simple simply connected groups defined and split over the prime field (such groups are also called {\it Chevalley groups}). Next, let $G = \mathrm{SO}_n(q)$ where $q$ is a
nondegenerate $n$-dimensional quadratic form over a field $K$ of characteristic $\neq 2$. Then $G$ is a 1-dimensional torus for $n = 2$, an absolutely almost simple group for $n = 3$ and $n \geq 5$, and a semi-simple but not absolutely almost simple group for $n = 4$. Furthermore, it is adjoint for odd $n$, and the corresponding simply connected group for all $n \geq 3$ is the spinor group $\mathrm{Spin}_n(q)$. It is $K$-split if and only if $q$ has the maximum Witt index, and is $K$-anisotropic (i.e., contains no $K$-split tori) if and only if $q$ does not represent zero over $K$.

\vskip1mm

Let us now describe the structure of an arbitrary connected algebraic group $G$.
First, $G$ possesses a maximal connected normal unipotent subgroup $R_u(G)$ called the {\it unipotent radical} of $G$. Then $R_u(H)$ is trivial for the quotient $H = G/R_u(G)$, i.e.
$H$ is \underline{{\it reductive}}.


Let now $G$ be reductive. Then its connected center $T = Z(G)^{\circ}$ is a torus, the derived subgroup $D = [G , G]$ is semi-simple, and the natural product map
$D \times T \to G$ is an isogeny. Moreover, if $G$ is defined over $K$, then $D$ and $T$ are also $K$-defined.

\section{Reduction of linear algebraic groups: Examples}\label{S-RedExamples}

In this section, we will discuss several concrete examples of reduction of algebraic groups defined over the field $\mathbb{Q}$ of rational numbers at rational primes. In these examples, we will just analyze the result of reduction of a set of defining equations with integer coefficients. These considerations will motivate the formal and more general definition of good reduction in the next section.



\vskip1mm

\noindent {\bf Example 3.1.} The defining equations for the $\mathbb{Q}$-groups $\mathrm{SL}_n$, $\mathrm{Sp}_{2n}$ and $\mathbb{D}_n$ from \S2.1 have integer coefficients, and reducing them modulo a prime $p$ we obtain the equations that define the groups of the same type over $\mathbb{F}_p$. Furthermore, let $G =
\mathrm{SO}_n(q)$ where $q = x_1^2 + \cdots + x_n^2$ $(n \geq 2)$. Then for any prime $p > 2$ the reduced equations define $\mathrm{SO}_n(\bar{q})$ where $\bar{q} = x_1^2 + \cdots + x_n^2$ over $\mathbb{F}_p$.

\vskip1mm

We note that all groups considered in this example are reductive, and so are their reductions. Here is an example of a different kind of behavior.

\vskip1mm

\noindent {\bf Example 3.2.} Fix a prime $p > 2$ and consider the binary quadratic form $q = x^2 - py^2$. Then $G = \mathrm{SO}_2(q)$ is a 1-dimensional $\mathbb{Q}$-anisotropic torus that consists of matrices of the form $\left( \begin{array}{cc} a & pb \\ b & a \end{array} \right)$ having determinant $1$.
Thus, the equations in terms of the matrix entries $x_{ij}$ that define $G$ are:
\begin{equation}\label{E:tor}
x_{11} = x_{22}, \ \ x_{12} = px_{21}, \ \ x_{11}^2 - px_{21}^2 = 1.
\end{equation}
Reducing these equations modulo $p$, we obtain the equations
$$
x_{11} = x_{22}, \ \  x_{12} = 0, \ \ x_{11}^2 = 1.
$$
The solutions to these equations are matrices of the form $\pm \left( \begin{array}{cc} 1 & 0 \\ \star & 1 \end{array}  \right)$. Thus, the reduced equations define an algebraic group over $\mathbb{F}_p$ which is disconnected and whose connected component is a 1-dimensional {\it unipotent} group! At the same time, reducing the equations (\ref{E:tor}) modulo any prime $q > 2$ different from $p$, we still get a 1-dimensional torus.

\vskip1mm

The group $G$ in this example can also be constructed as the norm  torus associated with the extension $L = \mathbb{Q}(\sqrt{p})$ of $\mathbb{Q}$ in the sense
that its $\mathbb{Q}$-points are precisely the elements of this extension having norm $1$. We will now explore the noncommutative version of this construction.

\vskip1mm

\noindent {\bf Example 3.3.} Let $p > 2$ be a prime and $D$ be the quaternion algebra over $\mathbb{Q}$ with basis $1, i, j, k$ and multiplication table
$$
i^2 = -1, \ j^2 = p, \ k = ij = -ji.
$$
We recall that the reduced norm of a quaternion $z = a + bi + cj + dk \in D$ is given by
$$
\mathrm{Nrd}_{D/\mathbb{Q}}(z) = a^2 + b^2 - pc^2 - pd^2.
$$
There exists an algebraic $\Q$-group $G$, usually denoted $\mathrm{SL}_{1 , D}$, whose $\mathbb{Q}$-points correspond to the quaternions $z \in D^{\times}$ with $\mathrm{Nrd}_{D/\mathbb{Q}}(z) = 1$. Using the regular representation of $D$ in the chosen basis, one can give an explicit matrix realization of $G$ and write down a set of defining equations with integer coefficients. One then proves that  (1) over $\overline{\Q}$, the group $\mathrm{SL}_{1,D}$ is isomorphic to $\mathrm{SL}_2$, hence is an absolutely almost simple group;  and (2) the reductions modulo $p$ of the defining equations yield a group that  has \emph{nontrivial} unipotent radical, i.e. is {\it not} reductive (see \cite[Example 3.5]{RR-Survey} for the details).

\vskip1mm

Thus, in contrast to Example 3.1, where the reductive groups had reductive reductions, in the last two examples, the reductive groups (in fact, a torus and an absolutely almost simple group) are defined by equations whose reductions yield non-reductive groups. This ``anomaly" is precisely what characterizes bad reduction! (We note that in a broad sense, the situation is analogous to the reduction of an elliptic curve given by (\ref{E:1}) at a prime dividing its discriminant in \S1, where the type of the curve changes.)

This discussion motivates the following conceptual approach to the notion of good reduction for reductive $\mathbb{Q}$-groups.
First, to make our notion ``local at $p$," instead of $\Z$, we take the localization $\mathbb{Z}_{(p)}$ of $\mathbb{Z}$ at the prime ideal $(p)$ as our base ring (concretely, $\mathbb{Z}_{(p)}$ is the subring of $\Q$ consisting of fractions in lowest terms whose denominators are prime to $p$).
Then we say that a reductive $\mathbb{Q}$-group has {\it good reduction} at a prime $p$ if it admits a matrix $\Q$-realization $G \hookrightarrow \mathrm{GL}_m$ where it can be defined by a system of polynomial equations with coefficients in $\mathbb{Z}_{(p)}$
such that the reduction of this system modulo $p$ still defines a reductive group; otherwise, we say that $G$ has {\it bad reduction} at $p$. Thus, the groups in Example 3.1  {\it do} have good reduction at the specified primes. On the other hand, Example 3.2 suggests that the 1-dimensional torus considered therein does not have good reduction at $p$ (which can, in fact, be established with some extra work), while it {\it does} have good reduction at all odd primes $q \neq p$.

For the group $G = \mathrm{SL}_{1 , D}$ considered in Example 3.3, the situation is more delicate. We have indicated that the realization of $G$ provided by the regular representation of $D$ always results in non-reductive reduction, and one can  show that $G$ indeed has
bad reduction at $p$ if $p \equiv 3(\mathrm{mod}\: 4)$. For $p \equiv 1(\mathrm{mod}\: 4)$, however, $G$ is $\Q$-isomorphic to $\mathrm{SL}_2$, hence has a realization that results in good reduction at $p$ (cf. the discussion in \S1 of why the elliptic curve $y^2 = x^3 - 625x$ actually has good reduction at $p = 5$). Thus, an adequate definition of good reduction should operate with the objects independent of the matrix realization of the given reductive group $G$ and the choice of defining equations. As we explained in \S 2, such an object is the Hopf algebra  $\mathbb{Q}[G]$ of $\Q$-regular functions, and the formal definition of good reduction in the next section is formulated precisely in terms of the existence of a $\mathbb{Z}_{(p)}$-structure on the latter whose reduction modulo $p$ represents a reductive group.


\section{Reductive algebraic groups with good reduction}\label{S-GoodRedDef}

\vskip2mm

Over general fields, one considers the notion of good reduction for reductive algebraic groups with respect to {\it discrete valuations}. So, we begin with a quick summary of the relevant definitions and some basic examples.

\vskip1mm

\subsection{Brief review of valuations.} First, we recall that a (normalized) discrete valuation on a field $K$ is a surjective map $v \colon K^{\times} \to \mathbb{Z}$ such that

\vskip1mm

(a) $v(ab) = v(a) + v(b)$;

\vskip1mm

(b) $v(a + b) \geqslant \min(v(a) , v(b))$ whenever $a + b \neq 0$.

\vskip1mm

\noindent One defines the {\it valuation ring} and the {\it valuation ideal}  respectively by $
\mathcal{O}(v) = \{ a \in K^{\times} \ \vert \ v(a) \geqslant 0 \} \cup \{ 0 \}$
and $
\mathfrak{p}(v) = \{ a \in K^{\times} \ \vert \ v(a) > 0 \} \cup \{ 0 \}.$
Then $\mathcal{O}(v)$ is a local ring with maximal ideal $\mathfrak{p}(v)$. The quotient ring $\mathcal{O}(v)/\mathfrak{p}(v)$ is called the {\it residue field} and will be denoted $K^{(v)}$. Furthermore, one considers the completion $K_v$ of $K$ with respect to the absolute value associated with $v$. Then $v$ naturally extends to a discrete valuation on $K_v$, which we will still denote by $v$. The corresponding valuation ring and valuation ideal in $K_v$ will be denoted $\mathcal{O}_v$ and $\mathfrak{p}_v$, respectively; we note that the quotient ring $\mathcal{O}_v/\mathfrak{p}_v$ coincides with the residue field $K^{(v)}$ defined above.

\vskip1mm

\noindent {\bf Example 4.1.}

\vskip1mm
\noindent (a) To every rational prime $p$, there corresponds the $p$-adic valuation on $\mathbb{Q}$ defined as follows: if $a \in \mathbb{Q}^{\times}$ is of the form $\displaystyle a = p^{\alpha} \cdot \frac{m}{n}$, with $m$ and $n$ relatively prime to $p$, then $v(a) = \alpha$. The corresponding valuation ring is the localization $\mathbb{Z}_{(p)}$, and the residue field is $\mathbb{F}_p$. Furthermore, the completion is the field of $p$-adic numbers $\mathbb{Q}_p$, and the valuation ring of $\mathbb{Q}_p$ is the ring $\mathbb{Z}_p$ of $p$-adic integers.

\vskip1mm

\noindent (b) Let $K = k(x)$ be the field of rational functions in one variable over a field $k$, and let $p(x) \in k[x]$ be a (monic) irreducible polynomial. Then the same construction as in part (a) enables us to associate to $p(x)$ a discrete valuation $v_{p(x)}$ on $K$. There is one additional discrete valuation on $K$ given by
$$
v_{\infty}\left(f/g\right) = \deg(g) - \deg(f), \ \ \mathrm{where} \ \ f , g \in k[x].
$$
Note that all of these valuations are trivial on the field of constants $k$, and cumulatively they constitute all valuations of $K$ with this property. These valuations are often called ``geometric" since they naturally correspond to the closed points of the projective line $\mathbb{P}^1_k$.

\vskip1mm

\noindent (c) Let again $K = k(x)$, but now assume that we are given a discrete valuation $v_0$ of $k$. Then $v_0$ can be extended to a discrete valuation $v$ on $K$ by first extending it to the polynomial ring $k[x]$ using the formula
$$
v(f(x)) := \min_{i = 0, \ldots , n} v_0(a_i)
$$
for $f(x) = a_n x^n + \cdots + a_0$
and then extending $v$ to $K$ by multiplicativity.

\subsection{Definition of good reduction.} We are finally ready to introduce a formal definition of good reduction for a reductive algebraic group $G$ defined over a field $K$ with respect to a discrete valuation $v$ of $K$. As we explained at the end of \S3, one needs to formulate it in terms of the Hopf algebra $K[G]$ of $K$-regular functions, i.e. exercising the point of view of the theory of affine group schemes. After stating the formal definition, however, we will exhibit a number of situations where good reduction can be characterized in very concrete terms.

In our discussion, we will use the following standard notation: given an affine scheme $X = \mathrm{Spec}\: A$, where $A$ is a commutative algebra over a commutative ring $R$, and a ring extension $R \subset R'$, we will denote by
$X \times_R R'$ the affine scheme $\mathrm{Spec}\: (A \otimes_R R')$ over $R'$ (usually referred to as ``base change").

\vskip1mm

\noindent {\bf Definition 4.2.} {\it Let $K$ be a field equipped with a discrete valuation $v$, and let $G$ be a connected reductive $K$-group. We say that $G$ has \emph{good reduction} at $v$ if there exists a reductive group scheme $\mathscr{G}$ over $\mathcal{O}_v$ with generic fiber $\mathscr{G} \times_{\mathcal{O}_v} K_v$ isomorphic to $G \times_K K_v$.}

\vskip1mm

Explicitly, this means that there should exist a Hopf $\mathcal{O}_v$-algebra $\mathcal{H}$ such that $\mathcal{H} \otimes_{\mathcal{O}_v} K_v \simeq K_v[G]$ (an ``$\mathcal{O}_v$-structure" on $K_v[G]$) and the algebra $\mathcal{H}/\mathfrak{p}_v\mathcal{H}$ over the residue field $K^{(v)}$ represents a connected reductive group (cf. the end of \S2.2).

We conclude this section with several examples of semi-simple groups where the conditions for good reduction can be described explicitly. These examples
will be sufficient for understanding the rest of the article.

\vskip1mm

\noindent {\bf Example 4.3.}

\noindent (a) It follows from the Chevalley construction that any absolutely almost simple simply connected $K$-{\it split} group $G$ has good reduction at every discrete valuation $v$ of $K$. (This generalizes the discussion of $\mathrm{SL}_n$ and $\mathrm{Sp}_{2n}$ in Example 2.1.)

\vskip1mm

\noindent (b) The construction of an absolutely almost simple group $G = \mathrm{SL}_{1 , D}$ associated with the group of norm 1 elements in a quaternion algebra $D$ can be generalized to any finite-dimensional central simple algebra $A$ over a field $K$ using the reduced norm homomorphism $\mathrm{Nrd}_{A/K} \colon A^{\times} \to K^{\times}$. We will discuss this group in more detail in Example 5.2(b) below, and now only indicate that $G = \mathrm{SL}_{1 , A}$ has good reduction at $v$ if and only if there exists an Azumaya $\mathcal{O}_v$-algebra $\mathcal{A}$ such that $A \otimes_K K_v \simeq \mathcal{A} \otimes_{\mathcal{O}_v} K_v$. This means that the algebra $A \otimes_K K_v$ should have an $\mathcal{O}_v$-structure $\mathcal{A}$ such that the quotient
$\mathcal{A}/\mathfrak{p}_v\mathcal{A}$ is a central simple algebra over $K^{(v)}$. Alternatively, the condition on $A$ can be expressed by saying that $A$ is {\it unramified} at $v$.

\vskip1mm

\noindent (c) Assuming that $\mathrm{char}\: K^{(v)} \neq 2$, the spinor group $G = \mathrm{Spin}_n(q)$ of a nondegenerate quadratic form $q$ in $n > 2$ variables over $K$ has good reduction at $v$ if and only if, over $K_v$, the form $q$ is equivalent to a quadratic form
$
\lambda (a_1 x_1^2 + \cdots + a_n x_n^2),
$
with
$\lambda \in K_v^{\times}$ and $a_i \in \mathcal{O}_v^{\times}$ for all $i = 1, \ldots , n.$

\section{Forms with good reduction and Galois cohomology}\label{S-Forms}

\subsection{$L/K$-forms.} The most general question in the spirit of Shafarevich's conjecture for reductive linear algebraic groups over arbitrary fields is the following:

\vskip1mm

\noindent {\it Let $K$ be a field equipped with a set $V$ of discrete valuations. What are the reductive algebraic groups of a given dimension that have good reduction at all $v \in V$? More specifically, what are the situations in which the number of $K$-isomorphism classes of such groups is finite?}

\vskip1mm

\noindent To make this question meaningful, one needs to specialize $K$ and $V$, which we will do in \S\ref{S-FinitenessConjecture}. However, we would first like to point out that considering reductive algebraic groups of a given dimension is far less natural than considering abelian varieties of the same dimension. Indeed, in a very coarse sense, all complex abelian varieties of dimension $d$ ``look the same": they are all analytically isomorphic to complex tori $\mathbb{C}^d / \Lambda,$ where $\Lambda \subset \mathbb{C}^d$ is a lattice of rank $2d$. At the same time, the ``fine" structure and classification of these varieties are highly involved as these varieties have nontrivial {\it moduli spaces}, which leads to infinite continuous families of nonisomorphic varieties.

On the contrary, reductive algebraic groups of the same dimension may look completely different; however, the structure theory (cf. \S2.2) implies
that over an algebraically (or separably) closed field, for each integer $d \geq 1$, there are only {\it finitely} many isomorphism classes of (connected) reductive groups of dimension $d$ (thus, there are no nontrivial moduli spaces for such groups). So, in the analysis of finiteness phenomena for reductive groups, it makes sense to focus on those classes of groups that becomes isomorphic over a separable closure of the base field; the main issue then becomes the passage from an isomorphism over the separable closure to an isomorphism over the base field (so-called {\it Galois descent}). This brings us to the following.

\vskip1mm

\noindent {\bf Definition 5.1.} {\it Let $G$ be a linear algebraic group over a field $K$, and let $L/K$ be a field extension. An algebraic $K$-group $G'$ is called an $L/K$-\emph{form} of $G$ if there exists an $L$-isomorphism $
G \times_K L \simeq G' \times_K L.$}

\vskip1mm

We will be interested mostly in $K^{\rm sep}/K$-forms of a linear algebraic group $G$, where $K^{\rm sep}$ is a fixed separable closure of $K$; these will often be referred to simply as {\it $K$-forms} of $G$. Here are some examples of forms of algebraic groups that we will encounter in the article.

\vskip1mm

\noindent {\bf Example 5.2.}

\vskip1mm

\noindent (a)  Let $T = \mathbb{D}_n$ be a $n$-dimensional $K$-split torus. Any other $n$-dimensional $K$-torus $T'$ splits over $K^{\mathrm{sep}}$, hence
$T \times_K K^{\mathrm{sep}} \simeq T' \times_K K^{\mathrm{sep}}$ over $K^{\mathrm{sep}}$.
This means that {\it all} $n$-dimensional $K$-tori are $K$-forms of $T$.

Similarly, let $G$ be an absolutely almost simple simply connected {\it split} algebraic group over $K$, and let $G'$ be any absolutely almost simple simply connected $K$-group of the same  type as $G$. Then $G'$ becomes split over $K^{\mathrm{sep}}$, entailing the existence of a $\Ks$-isomorphism
$G \times_K K^{\mathrm{sep}} \simeq G' \times_K K^{\mathrm{sep}}$.
Thus, $G'$ is a $K$-form of $G$, and hence every absolutely almost simple simply connected $K$-group is a $K$-form of an absolutely almost simple simply connected split $K$-group.

\vskip1mm

\noindent (b) Let $A$ be a central simple $K$-algebra of degree $n$, and let $G = \mathrm{SL}_{1 , A}$ be the algebraic $K$-group associated with the group of elements of reduced norm 1 in $A$. There is an isomorphism of $K^{\mathrm{sep}}$-algebras
$
\varphi \colon A \otimes_K K^{\mathrm{sep}} \longrightarrow \mathrm{M}_n(K^{\mathrm{sep}}),
$
and in terms of this isomorphism, the reduced norm of $a \in A$ is given by $\mathrm{Nrd}_{A/K}(a) = \det(\varphi(a \otimes 1))$. It follows that the group $G \times_K K^{\mathrm{sep}}$ is $K^{\mathrm{sep}}$-isomorphic to $\mathrm{SL}_n$; in other words, $G$ is a $K$-form of $\mathrm{SL}_n$. In fact, it is an {\it inner form}, and all inner forms are obtained this way --- see \cite[Ch. II, \S2.3.4]{Pl-R} for the details.

\vskip1mm

\noindent (c) Assume that $\mathrm{char}\: K \neq 2$. Let $q$ be a nondegenerate quadratic form in $n \geq 3$ variables over $K$, and let
$G = \mathrm{Spin}_n(q)$ be the corresponding spinor group.
Then any other nondegenerate quadratic form $q'$ over $K$ in $n$ variables is equivalent to $q$ over $K^{\mathrm{sep}}$, implying that $G' = \mathrm{Spin}_n(q')$ is a $K$-form of $G$. If $n$ is odd, then the groups $G' = \mathrm{Spin}_n(q')$ account for all $K$-forms of $G$. However, for $n$ even, there may be other $K$-forms of $G$ defined in terms of the (universal covers of) unitary groups of (skew-)hermitian forms over noncommutative central division $K$-algebras with an involution of the first kind (such as, for example, quaternion algebras).

\vskip1mm

We can now give a more precise version of the above question concerning reductive groups with good reduction.

\vskip1mm

\noindent {\bf Question 5.3.} {\it  Let $K$ be a field equipped with a set $V$ of discrete valuations. What are the $K$-forms (or even just inner $K$-forms --- see Example 5.4 below) of a given reductive algebraic $K$-group $G$ that have good reduction at all $v \in V$? More specifically, in what situations is the number of $K$-isomorphism classes of such $K$-forms \emph{finite}?}

\vskip1mm

We will return to this question in the next section and formulate precise conjectures in several situations. First, however, we would like to review
how the forms of algebraic groups (and other objects) are classified in terms of Galois cohomology.


\subsection{Classification of forms and Galois cohomology.}\label{S-GaloisCoh}

Let $\Gamma$ be a group, and $A$ be a $\Gamma$-module, i.e. an abelian group equipped with a $\Gamma$-action. In homological algebra, one defines the cohomology groups $H^i(\Gamma , A)$ for all $i \geq 0$. This definition can be extended to a not necessarily abelian $\Gamma$-group $A$ for $i = 0, 1$, in which case $H^1(\Gamma , A)$ does not have a natural group structure and is treated as a set with a distinguished element -- the equivalence class of the trivial cocycle. We refer the reader to \cite{Serre-GC} for the description of this construction and its basic properties, and will focus here on the applications to the classification of forms.

Given a finite Galois extension $L/K$ with Galois group $\Gamma$, there are two sources of $\Gamma$-groups. First, for any algebraic $K$-group $G$, the group $\Gamma$ naturally acts on the group of $L$-points $G(L)$. In this case, instead of $H^1(\Gamma , G(L))$ we write $H^1(L/K , G)$. Second, we can look at the group $A_L$ of $L$-defined automorphisms of $G$. Then again $A_L$ receives a natural action of $\Gamma$ and we may consider the set $H^1(\Gamma , A_L) = H^1(L/K , A_L)$. All this generalizes to infinite Galois extensions, but then one should regard the Galois group $\Gamma$ as a profinite group and define $H^1$ using continuous 1-cocycles on $\Gamma$ with values in a discrete group. The most important case here is when $L = K^{\mathrm{sep}}$ so that $\Gamma$ is the absolute Galois group of $K$. We then abbreviate $H^1(K^{\mathrm{sep}}/K , G)$ to just $H^1(K , G)$.

Let now $L/K$ be an arbitrary Galois extension (finite or infinite). For an algebraic $K$-group $G$, we let $F(L/K , G)$ denote the set of $K$-isomorphism classes of $L/K$-forms of $G$. One shows that there is a natural bijection $\gamma_{L/K} \colon F(L/K , G) \to H^1(L/K , A_L)$; in particular, we have a natural bijection $\gamma_K \colon F(K^{\mathrm{sep}}/K , G) \to H^1(K , A_{K^{\mathrm{sep}}})$. The examples below demonstrate how this cohomological description enables one to obtain a classification of forms in more explicit terms.

\vskip1mm

\noindent {\bf Example 5.4.}

\vskip.2mm

\noindent (a) Let $T = \mathbb{D}_n$ be the $n$-dimensional $K$-split torus. Set $\Gamma = \mathrm{Gal}(K^{\mathrm{sep}}/K)$. Then the automorphism group $A_{K^{\mathrm{sep}}}$ as $\Gamma$-group can be identified with $\mathrm{GL}_n(\mathbb{Z})$ equipped with the trivial action. Since any $n$-dimensional $K$-torus is a $K$-form of $T$ (cf. Example 5.2(a)), it follows that the $K$-isomorphism classes of $n$-dimensional $K$-tori are in bijection with the conjugacy classes
of continuous homomorphisms $f \colon \Gamma \to \mathrm{GL}_n(\mathbb{Z})$.

\vskip1mm

\noindent (b) Let $G$ be an absolutely almost simple simply connected $K$-split group. Then the $K$-forms of $G$ are precisely all absolutely almost simple simply connected $K$-groups of the same type as $G$, so their $K$-isomorphism classes are classified by the elements of $H^1(K , A_{K^{\mathrm{sep}}})$. In the case at hand, the automorphism group $A_{K^{\mathrm{sep}}}$ is the group of $K^{\mathrm{sep}}$-points of an algebraic $K$-group $A$, which is an extension of the group of inner automorphisms (typically identified with the adjoint group $\overline{G}$) by the finite group of all symmetries of the Dynkin diagram of $G$. Then the $K$-forms that correspond to the cohomology classes lying in the image of the natural map
$
H^1(K , \overline{G}) \longrightarrow H^1(K , A_{K^{\mathrm{sep}}})
$
are called {\it inner forms}; all other forms are called {\it outer forms}.

\vskip1mm

A similar approach applies to the classification of $K$-forms of other ``objects," which provides an interpretation of $H^1$ in certain situations.

\noindent {\bf Example 5.5.}

\vskip.2mm

\noindent (a) Let $D = \mathrm{M}_n(K)$ be the matrix algebra of degree $n$ over $K$. Then any other central simple algebra $D'$ of
degree $n$ over $K$ is a $K$-form of $D$ in the obvious sense. It follows from the Skolem-Noether theorem that for any field $F$, the group of $F$-automorphisms
of $\mathrm{M}_n(F)$ is $\mathrm{PGL}_n(F)$. Thus, we conclude that the isomorphism classes of central simple $K$-algebras of degree $n$ are in bijection with the elements of $H^1(K , \mathrm{PGL}_n)$.

\vskip1mm

\noindent (b) Fix a nondegenerate quadratic form $q$ in $n$ variables over a field $K$ of characteristic $\neq 2$. Then any other nondegenerate quadratic form $q'$ in $n$ variables over $K$ is a $K$-form of $q$. So, the equivalence classes of nondegenerate quadratic forms in $n$ variables over a field $K$ are in bijection with the elements of $H^1(K , \mathrm{O}_n(q))$ for the corresponding orthogonal group.

\section{The finiteness conjecture for reductive groups with good reduction}\label{S-FinitenessConjecture}

We now return to Question 5.3 about the forms of a reductive group having good reduction at a given set of discrete valuations and discuss two natural choices for the field $K$ and the set of discrete valuations $V$ of $K$.

\subsection{The Dedekind case.}  Here, the field $K$ is the fraction field of a Dedekind ring $R$ and $V$ is the set of discrete valuations of $K$ associated with the maximal ideals of $R$. Early interest in this case can be traced back to the work of G.~Harder (1967)  and J.-L.~Colliot-Th\'el\`ene and  J.-J.~Sansuc (1979). The basic case where $K = \mathbb{Q}$ and $R = \mathbb{Z}$, and hence $V$ is the set of all $p$-adic valuations of $\Q$, was considered by B.H.~Gross \cite{Gross}  and B.~Conrad (2014). In particular, Gross showed that an absolutely almost simple simply connected $\Q$-group $G$ that has a good reduction at all primes $p$ necessarily splits over all $\mathbb{Q}_p$. Furthermore, the following finiteness result is true over any number field.
\begin{prop}\label{P:XXX1}
Let $G$ be an absolutely almost simple simply connected algebraic group over a number field $K$, and let $V^K_f$ denote the set of all discrete valuations
of $K$. Then for any subset $V \subset V_f^K$  such that $V_f^K \setminus V$ is finite, the number of $K$-isomorphism classes of $K$-forms of $G$ that have good reduction at all $v \in V$ is finite.
\end{prop}

Another basic example of a Dedekind ring is the polynomial ring $R = k[x]$ over a field $k$. We let $V$ denote the set of discrete valuations of the field of rational functions $K = k(x)$ associated with monic irreducible polynomials $p(x) \in k[x]$ (see Example 4.1(b)). We then have the following result, which was first discovered by M.S.~Raghunathan and A.~Ramanathan (1984).
\begin{thm}\label{T:XXX2}
Let $G_0$ be a (connected) semi-simple simply connected algebraic group over a field $k$. If $G'$ is a $K$-form of the group $G = G_0 \times_k K$ that has good reduction at all $v \in V$ and splits over $k^{\rm sep}(x)$, then $G' = G'_0 \times_k K$ for some $k$-form $G'_0$ of $G_0$.
\end{thm}

Yet another case where forms with good reduction have been considered in full involves the ring of Laurent polynomials $R = k[x , x^{-1}]$. Here, the set $V$ consists of the discrete valuations of $K = k(x)$ corresponding to monic irreducible polynomials $p(x) \in k[x]$, with $p(x) \neq x$. In this case, a description of the $K$-forms of the group $G$ as in Theorem \ref{T:XXX2} that have good reduction at all $v \in V$ was
given by V.I.~Chernousov, P.~Gille, and A.~Pianzola (2012), under some minor restrictions on the characteristic.
This result is too technical to be stated here, so we only indicate that it played a crucial role in the proof of the conjugacy of the analogues of Cartan subalgebras in certain infinite-dimensional Lie algebras (we refer the reader to \cite[\S5.1]{RR-Survey} for precise statements and references).

Generalizing the last two examples, one can consider the (Dedekind) ring $R = k[C]$ of regular functions of a smooth geometrically integral affine curve $C$ over a field $k$, in which case $V$ consists of the discrete valuations of the function field $K = k(C)$ associated with the closed points of $C$. It seems that no explicit description of the $K$-forms of a given reductive $K$-group $G$ that have good reduction at all $v \in V$ is available in most situations. So, we will instead focus on the qualitative question concerning the conditions that ensure the finiteness of the number of $K$-isomorphism classes of such forms. We observe that if $G_0$ is a reductive algebraic $k$-group and $G = G_0 \times_k K$ is the base change of $G_0$ to $K$, then for any $k$-form $G'_0$ of $G_0$, the group $G' = G'_0 \times_k K$ is a $K$-form of $G$ that has good reduction at all $v \in V$. So, even though non-isomorphic $k$-forms may become isomorphic after base change to $K$, in order to have an affirmative answer to our question in a sufficiently general situation, one needs to assume that over $k$, the groups at hand have only finitely many non-isomorphic forms. This basically amounts to a hypothesis on the finiteness of Galois cohomology.

In \cite[Ch. III]{Serre-GC}, Serre described a class of fields over which one does have the required finiteness. For this he first
introduced the following {\it condition} (F) on a profinite group $\mathscr{G}$:

\vskip1mm

\noindent {\it For every integer $m \ge 1$, the group $\mathscr{G}$ has only finitely many open subgroups of index $m$.}

\vskip1mm

\noindent He then defined a field $K$ to be {\it of type} (F) if it is {\it perfect} and its absolute Galois group $\mathrm{Gal}(K^{\mathrm{sep}}/K)$ satisfies (F). Examples of such fields include $\R$ and $\C$, finite fields $\mathbb{F}_q$, and $p$-adic fields $\Q_p$, among many others. The key result is that if $K$ is a field of type (F), then the set $H^1(K,G)$ is finite for any linear algebraic $K$-group $G$.  Moreover, if $K$ is of characteristic 0 and of type (F), then any linear algebraic group has finitely many $K$-forms. Based on these results, we would like to propose the following finiteness conjecture for forms with good reduction.

\vskip1mm

\noindent {\bf Conjecture 6.3.} {\it Let $k$ be a field of characteristic 0 and of type $(\mathrm{F})$. Suppose $K = k(C)$ is the function field of a smooth geometrically integral affine curve $C$ over $k$, and let $V$ be the set of discrete valuations of $K$ associated with the closed points of $C$.
Then for an absolutely almost simple simply connected algebraic $K$-group $G$, the number of $K$-isomorphism classes of $K$-forms $G'$ of $G$ that have good reduction at all $v \in V$ is finite.}

\vskip1mm

We refer the reader to \cite[\S5.2]{RR-Survey} for a more general form of this conjecture.


\subsection{The finiteness conjecture for arbitrary finitely generated fields.} We now turn to the main finiteness conjecture for forms with good reduction over arbitrary finitely generated fields. Let us point out that, unlike the Dedekind situation discussed above, forms with good reduction had not been considered in the higher-dimensional setting prior to our work.

We first note that an arbitrary finitely generated field $K$ possesses natural sets of discrete valuations called {\it divisorial}. In geometric language, such sets consist of the discrete valuations associated with prime divisors on a {\it model} $\mathfrak{X}$ of $K$, i.e. a normal separated irreducible scheme of finite type over $\Z$ (if $\mathrm{char}\: K = 0$) or over a finite field (if $\mathrm{char}\: K > 0$) such that $K$ is the field of rational functions on $\mathfrak{X}$. Algebraically, this construction amounts to finding a subring $R$ of $K$ whose fraction field is $K$ and such that $R$ is integrally closed (in $K$) and a finitely generated $\mathbb{Z}$-algebra (or $\mathbb{F}_p$-algebra). For any height one prime ideal $\mathfrak{p} \subset R$ (i.e. a minimal nonzero prime ideal), the localization $R_{\mathfrak{p}}$ is a discrete valuation ring, hence gives rise to a discrete valuation $v_{\mathfrak{p}}$ of $K$. Then $V$ consists of such valuations $v_{\mathfrak{p}}$ corresponding to all height one prime ideals of $R$. Furthermore, we observe that any two divisorial sets $V_1$ and $V_2$ associated with two different models of $K$ are {\it commensurable}, i.e. $V_i \setminus (V_1 \cap V_2)$ is finite for $i = 1, 2$ (which makes a divisorial set of places almost canonical), and that for any finite subset $S$ of a divisorial set $V$, the set $V \setminus S$ contains a divisorial set. Here is an example that gives a concrete illustration of these notions.

\vskip1mm

\noindent {\bf Example 6.4.} The field of rational functions $K = \mathbb{Q}(x)$ is the fraction field of $R = \mathbb{Z}[x]$.
It is well-known that any height one prime ideal $\mathfrak{p}$ of $R$ is  principal with a generator $p$ of one of the following two types: (a) $p \in \mathbb{Z}[x]$ is an irreducible polynomial with content 1; or (b) $p \in \mathbb{Z}$ is a rational prime. In the first case,
the associated discrete valuation of $K$ coincides with the one corresponding to the monic irreducible polynomial in $\mathbb{Q}[x]$ obtained by dividing $p$ by its leading coefficient (cf. Example 4.1(b)). We will call such discrete valuations {\it geometric} and denote the set of all of these valuations by $V_0$. In the second case, the associated discrete valuation coincides with the extension of the $p$-adic valuation described in Example 4.1(c). We will refer to such valuations as {\it arithmetic} and denote the set of all of these valuations by $V_1$. Thus, the set of divisorial valuations of $K$ associated with the model $\mathfrak{X} = \mathrm{Spec}\: R$ is $V = V_0 \cup V_1$.

\vskip1mm

We now come to our central conjecture for forms with good reduction over arbitrary finitely generated fields.

\vskip1mm

\noindent {\bf Conjecture 6.5. (Main Conjecture for forms with good reduction)} {\it Let $G$ be a (connected) reductive algebraic group over a finitely generated field $K$, and let $V$ be a divisorial set of discrete valuations of $K$. Then the set of $K$-isomorphism classes of (inner) $K$-forms $G'$ of $G$ that have good reduction at all $v \in V$ is finite (at least when the characteristic of $K$ is ``good'').}

\vskip1mm

As follows from our discussion, in a broad sense this conjecture should be viewed as an analogue of Shafarevich's conjecture in the context of reductive algebraic groups over higher-dimensional fields. It appears to play a key role in the analysis of other finiteness properties in this situation, and in \S\ref{S-Hasse} we will discuss its connections with the local-global principle. Due to space limitations, we will not be able to presents its links to the {\it genus problem} (cf. \cite[\S8]{RR-Survey}), but in \S9 we will describe an application of this work to length-commensurable Riemann surfaces.

\section{Finiteness aspects of local-global properties}\label{S-Hasse}

We now discuss connections between Conjectures 6.3 and 6.5 and local-global principles.
Historically, one of the earliest local-global results was the Hasse-Minkowski theorem, which states that given a quadratic form $q(x_1, \dots, x_n)$ over a number field $K$, the equation
$
q(x_1, \dots, x_n) = 0
$
has a nonzero solution in $K^n$ if and only if it has a nonzero solution in $K_v^n$ for all valuations $v$ of $K$ (including the archimedean ones).
The essence of this result is well encapsulated in the general idea of the local-global principle, which can be formulated as follows:

\vskip1mm

\noindent {\it Suppose $K$ is a field equipped with a set $V$ of (rank one) valuations, and let $\mathcal{C}$ be a class of algebraic $K$-varieties. One says that \emph{the local-global principle holds for $\mathcal{C}$ with respect to $V$} if any variety $X \in \mathcal{C}$ has a rational point over $K$ if and only if it has a rational point over the completion $K_v$ for every $v \in V$.}

\vskip1mm

\noindent From this perspective, the Hasse-Minkowski theorem simply states that the class of projective quadrics over any number field $K$ satisfies the local-global principle with the respect to the set $V$ of all valuations of $K$.

In this article, we will focus on the cohomological approach to the local-global principle, where one considers an algebraic group $G$ defined over a field $K$ equipped with a set $V$ of valuations and then analyzes the global-to-local map in Galois cohomology $
\theta_{G,V} \colon H^1(K,G) \to \prod_{v \in V} H^1(K_v, G).$ One then defines the {\it Tate-Shafarevich set} $\Sha(G,V)$ as
the kernel $\ker \theta_{G,V}$, i.e. the preimage of the distinguished element. It is easy to see that the validity of the local-global principle for the class of principal homogeneous spaces of $G$ (or $G$-torsors) is equivalent to the triviality of $\Sha(G,V)$, but in fact one is even more interested in the injectivity of $\theta_{G,V}$ which for a noncommutative $G$ is {\it not} an automatic consequence of this property. The reason is that the injectivity enables one to give a local-global treatment of various {\it equivalences}. For example, it follows from Example 5.5(b) that the injectivity of $\theta_{G,V}$ for the orthogonal group $G = \mathrm{O}_n(q)$ is equivalent to the well-known consequence of the Hasse-Minkowski theorem is that two nondegenerate quadratic forms $q_1(x_1, \dots, x_n)$ and $q_2(x_1, \dots, x_n)$ over a number field $K$ are equivalent over $K$ if and only if they are equivalent over the completions $K_v$ for all valuations $v$ of $K$ (including the archimedean ones). So, one says that  the {\it cohomological local-global principle} holds for $G$ with respect to $V$ if the  map $\theta_{G,V}$ is injective. Here is what this amounts to in the classical situation.

\vskip1mm

\noindent {\bf Example 7.1.}


\noindent (a) The cohomological principle for the norm torus $T$ associated with the group of norm 1 elements in a finite separable extension $L/K$ is equivalent to the statement that  an element $a \in K^{\times}$ is a norm in the extension $L/K$  if and only if it is a norm locally at all $v \in V$. According to the famous Hasse Norm Theorem, this is indeed the case for all cyclic Galois extensions of number fields when $V$ is the set of all valuations of $K$, including the archimedean ones.

\vskip1mm

\noindent (b) Let $G = {\rm PGL}_n.$ It follows from Example 5.5(a) that the cohomological local-global principle for $G$ amounts to the fact
that two central simple algebras $A_1$ and $A_2$ of degree $n$ over $K$ are isomorphic if and only the algebras $A_1 \otimes_K K_v$ and $A_2 \otimes_K K_v$ are isomorphic over $K_v$ for all $v \in V$. The truth of this for all $n \geq 2$ makes the natural map of Brauer groups
$\mathrm{Br}(K) \longrightarrow \prod_{v \in V} \mathrm{Br}(K_v)$ injective, which is indeed the case if  $K$ is a number field and $V$ is the set of all valuations of $K$ including the archimedean ones according to the Albert-Hasse-Brauer-Noether theorem.

\vskip1mm

Starting with these examples,  the study of the cohomological local-global principle was eventually extended to arbitrary linear algebraic groups. One of the crowning achievements of the arithmetic theory of algebraic groups was the proof of the fact that the cohomological local-global principle holds for \emph{all} semi-simple simply connected and adjoint groups over number fields (cf. \cite[Ch. 6]{Pl-R}).
While the cohomological local-global principle fails for general semi-simple groups (see \cite[Ch. III]{Serre-GC}), it was nevertheless proved by Borel and Serre using reduction theory for adelic groups that when $K$ is a number field, the map $\theta_{G,V}$ is \emph{proper} (i.e., the pre-image of a finite set is finite) for \emph{any} linear algebraic group $G$ whenever $V$ contains almost all valuations of $K$.
Informally, these results mean that while the local-global principle for groups over number fields may fail, the deviation from it is always of finite size.
It should be noted that the finiteness of the Tate-Shafarevich group for abelian varieties over number fields remains one of the major open problems. At the same time, for linear algebraic groups one expects such finiteness properties to hold much more generally.

\vskip1mm

\noindent {\bf Conjecture 7.2.} {\it Let $G$ be a reductive algebraic group defined over a finitely generated field $K$, and let $V$ be a divisorial set of discrete valuations of $K$. Then the global-to-local map $\theta_{G , V}$ is \emph{proper}. In particular, the Tate-Shafarevich set $\Sha(G , V)$ is finite.}

\vskip1mm

We note that it also makes sense to formulate an analogous conjecture in the case where $k$ is a field of characteristic 0 and of type (F) (cf. \S6.1), $K= k(C)$ is the function field of a smooth geometrically integral affine curve $C$ over $k$, and $V$ is the set of discrete valuations of $K$ associated with the closed points of $C$.

We will discuss some of the available results on these conjectures in the next section, but to conclude this section, let us mention the following rather surprising connection between the analysis of groups with good reduction and local-global principles: The truth of Conjecture 6.5 for the inner forms of an absolutely almost simple simply connected group $G$ automatically implies the properness of the global-to-local map
$
\theta_{\overline{G} , V} \colon H^1(K , \overline{G}) \longrightarrow \prod_{v \in V} H^1(K_v , \overline{G})
$
for the corresponding adjoint group $\overline{G}$ (cf. \cite{R-ICM}). Thus, Conjectures 6.5 and 7.2 are related.

\section{Overview of results}\label{S-Results}

We will now give a brief overview of some available results on Conjectures 6.3, 6.5, and 7.2, referring the interested reader to \cite{RR-Survey} for a more detailed discussion.


\vskip2mm

\subsection{Algebraic tori.} Recently, all these conjectures were resolved for algebraic tori \cite{RR-Tori}. First, it was shown that if $K$ is a finitely generated field of characteristic zero and $V$ is a divisorial set of places of $K$, the for any  $n \geq 1$, the set of $K$-isomorphism classes of $n$-dimensional $K$-tori that have good reduction at all $v \in V$ is finite. One should view this result that settles Conjecture 6.5 for tori as a far-reaching generalization of the Hermite-Minkowski theorem that a given number field has only finitely many extensions of a given degree that ramify only at a given finite set of places. The proof relies on the description of the isomorphism classes of $n$-dimensional tori given in Example 5.4(a) in conjunction with the fact  that the \'etale fundamental group of the model $\mathfrak{X}$ of $K$ used to define $V$ satisfies Serre's condition (F).
A similar result is also valid in the context of Conjecture 6.3: If $K = k(X)$ is the function field of a smooth irreducible algebraic variety $X$ over a field $k$ that has characteristic zero and is of type (F) and $V$ is the set of discrete valuations of $K$ associated with the prime divisors of $X$, then for any $n \geq 1$ there are only finitely many $K$-isomorphism classes of $n$-dimensional $K$-tori that have good reduction at all $v \in V$.

We now turn to the local-global principles for tori, noting that since tori are commutative, the properness of the global-to-local map for them is equivalent to the finiteness of the Tate-Shafarevich group. The following theorem confirms Conjecture 7.2 for tori.
\begin{thm}\label{T:Res2}
Let $K$ be a finitely generated field and $V$ be a divisorial set of discrete valuations of $K$. Then for any algebraic $K$-torus $T$, the Tate-Shafarevich group
$$
\Sha(T , V) = \ker \left( H^1(K , T) \to \prod_{v \in V} H^1(K_v , T) \right)
$$
is finite.
\end{thm}

We recall that the classical proof of the finiteness of $\Sha(T , V)$ when $K$ is a number field and $V$ is set of {\it all} places of $K$ relies on Tate-Nakayama duality, which is not available in the general situation. One of the proofs of Theorem \ref{T:Res2} in \cite{RR-Tori} systematically uses adelic techniques in the context of arbitrary finitely generated fields and their divisorial sets of valuations, which, to the best of our knowledge, have not been previously employed in this generality. In particular, it demonstrates that when $K$ is a number field, the finiteness of the Tate-Shafarevich group can be established without Tate-Nakayama duality, and is actually a consequence of two basic facts: the finite generation of the group of $S$-units and the finiteness of the class number. The other proof requires an additional assumption on the characteristic of $K$, but applies also to some fields that are not finitely generated, such as function fields of algebraic varieties over fields of characteristic zero and having type (F).

\vskip2mm

\subsection{Absolutely almost simple groups.}

First, Conjecture 6.5. has been proven completely for inner forms of type $\mathsf{A}_n.$
\begin{thm}\label{T:ResX1}
Let $K$ be a finitely generated field, $V$ a divisorial set of discrete valuations of $K$, and $n \geq 2$ an integer prime to $\mathrm{char}\: K$. Then the number of $K$-isomorphism classes of groups of the form $\mathrm{SL}_{1,A}$, where $A$ is a central simple $K$-algebra of degree $n$, that have good reduction at all $v \in V$, is finite.
\end{thm}


Next, we have resolved Conjecture 6.5 for groups of several other types over some special classes of finitely generated fields, particularly over {\it two-dimensional global fields} \cite{CRR-Spinor}. Following Kato, by a two-dimensional global field, we mean either the function field $K = k(C)$ of a smooth geometrically integral curve $C$ over a number field $k$, or the function field $K = k(S)$ of a smooth geometrically integral surface $S$ over a finite field $k$.

\begin{thm}\label{T:ResX2}
Let $K$ be a two-dimensional global field of characteristic $\neq 2$ and let $V$ be a divisorial set of discrete valuations of $K$. Fix an integer $n \geq 5$. Then the set of $K$-isomorphism classes of spinor groups $G = \mathrm{Spin}_n(q)$ of nondegenerate quadratic forms in $n$ variables over $K$ that have good reduction at all $v \in V$ is finite.
\end{thm}

The proof uses a variety of techniques, including Milnor's conjecture on quadratic forms (proved by Voevodsky). Similar statements are also available for some special unitary groups of types $\textsf{A}_n$ and $\textsf{C}_n$ and for groups of type $\textsf{G}_2$. Without going into the details, we should  also point out that the considerations involved in the proof of Theorem \ref{T:ResX2} lead to analogous results concerning
Conjecture 6.3.

As we noted at the end of \S\ref{S-Hasse}, Conjecture 6.5 automatically implies Conjecture 7.2 for adjoint groups. Thus, it follows from Theorem \ref{T:ResX1} that if $K$ is a finitely generated field equipped with a divisorial set of discrete valuations $V$, the global-to-local map $\theta_{G , V}$ in Galois cohomology is proper for $G = \mathrm{PSL}_{1 , A}$ for any central simple algebra $A$ of degree $n$ (thus, in particular, for $G = \mathrm{PSL}_n$) when $n$ is prime to the characteristic. Furthermore, for a two-dimensional global field $K$ of characteristic $\neq 2$ and its divisorial set of places $V$, the properness of $\theta_{G , V}$ is known in the following cases:

\vskip1mm

\noindent $\bullet$ \parbox[t]{16.5cm}{$G = \mathrm{SL}_{1,A}$, where $n$ is a square-free integer prime to $\mathrm{char}\: K$ and $A$ is a central simple $K$-algebra of degree $n$;}

\vskip1mm

\noindent $\bullet$ \parbox[t]{16.5cm}{$G = \mathrm{SO}_n(q)$, where $q$ is a nondegenerate quadratic form over $K$ in $n \geq 5$ variables;}

\vskip1mm

\noindent $\bullet$ \parbox[t]{16.5cm}{$G$ is a simple group of type $\mathsf{G}_2$}

\vskip2mm

\noindent (Note that for ${\rm SO}_n(q)$, the above statement follows from Theorem \ref{T:ResX2} only for {\it odd} $n$.)

Our recent analysis of the unramified cohomology of function fields
of rational surfaces and certain Severi-Brauer varieties over number fields has yielded
similar properness statements in these cases as well (cf. \cite{RR-Tori}). The investigation of $\theta_{G , V}$ for other groups over arbitrary finitely
generated fields continues.

\section{An application to Riemann surfaces}\label{S-RiemannSurfaces}

Let $\mathbb{H}$ be the complex upper half-plane equipped with the standard hyperbolic metric  and the action of $\mathrm{SL}_2(\mathbb{R})$ by fractional linear transformations. We are interested in Riemann surfaces of the form $M = \Gamma \backslash \mathbb{H}$, where $\Gamma \subset \mathrm{SL}_2(\mathbb{R})$ is a discrete subgroup with torsion-free image in $\mathrm{PSL}_2(\mathbb{R})$ --- recall that all compact Riemann surfaces of genus $g > 1$ admit such a description. If $\Gamma$ is an arithmetic Fuchsian group, then $M$ is said to be {\it arithmetically defined}. The investigation of Riemann surfaces proceeds in different directions, one of which ({\it spectral geometry}) focuses on characterizing Riemann surfaces (and other types of spaces) in terms of various {\it spectra}. A classical example
is the Laplace spectrum $E(M)$,
leading to M.~Kac's famous question {\it Can one hear the shape of a drum}? It is also natural
to consider the {\it length spectrum} $L(M)$, defined as the collection of the lengths of all closed geodesics in $M$.

While neither spectrum determines a Riemann surface $M$ up to isometry (cf. \cite{Sunada}, \cite{Vigneras}), A.~Reid \cite{Reid} established the following commensurability statement: for two {\it arithmetically defined} Riemann surfaces $M_1$ and $M_2$, the condition of iso-length-spectrality $L(M_1) = L(M_2)$ implies that $M_1$ and $M_2$ are {\it commensurable}, i.e. have a common finite-sheeted cover.
In fact, in the arithmetic situation, commensurability  also follows from the much weaker condition $\mathbb{Q} \cdot L(M_1) = \mathbb{Q} \cdot L(M_2)$, which is referred to as {\it length-commensurability}.

To establish commensurability, Reid considers the algebra $A_{\Gamma}$ attached to a Riemann surface $M = \Gamma \backslash \mathbb{H}$, defined as the  $\mathbb{Q}$-subalgebra of the matrix algebra $\mathrm{M}_2(\mathbb{R})$ spanned by $\Gamma^{(2)}$, the subgroup of $\Gamma$ generated by the squares of all elements. It turns out that for any Zariski-dense subgroup $\Gamma \subset \mathrm{SL}_2(\mathbb{R})$, $A_{\Gamma}$ is a quaternion algebra,
whose center is the {\it trace field} $K_{\Gamma} := \mathbb{Q}(\mathrm{tr}\: \gamma \, \vert \, \gamma \in \Gamma^{(2)})$.  If $\Gamma_1$ and $\Gamma_2$ are commensurable, then $A_{\Gamma_1} = A_{\Gamma_2}$,
making $A_{\Gamma}$ an invariant of the commensurability class of $\Gamma$. Moreover, for $\Gamma$ an arithmetic group, $A_{\Gamma}$ actually determines $\Gamma$ up to commensurability.
Reid first shows that two iso-length-spectral arithmetically defined Riemann surfaces $M_1 = \Gamma_1 \backslash \mathbb{H}$ and $M_2 = \Gamma_2 \backslash \mathbb{H}$ have the same trace field $K_{\Gamma_1} = K_{\Gamma_2} =: K$, which is a number field. Then he uses the Albert-Hasse-Brauer-Noether Theorem (AHBN) to conclude that $A_{\Gamma_1} \simeq A_{\Gamma_2}$ over $K$, which yields the required result.

``Most" Riemann surfaces, however, are {\it not} arithmetic, and their investigation presents many challenges. In our set-up, while $A_{\Gamma}$ remains one of very few available invariants of the commensurability class of $\Gamma$, it no longer determines $\Gamma$ up to commensurability. What one can prove
assuming only that the subgroups $\Gamma_1, \Gamma_2 \subset \mathrm{SL}_2(\mathbb{R})$ are Zariski-dense, is that the length-commensurability of $M_1 = \Gamma_1 \backslash \mathbb{H}$ and $M_2 = \Gamma_2 \backslash \mathbb{H}$ implies
that $K_{\Gamma_1} = K_{\Gamma_2} = K$. In contrast to the arithmetic case, however,  this field may be an arbitrary finitely generated subfield of $\mathbb{R}$, where no analogues of (ABHN) are known or can even be established. Instead, using good reduction and Theorem \ref{T:ResX1}
we were able to prove the following statement that does not involve any arithmeticity or even finite volume assumptions.
\begin{thm}\label{T:RSurf}
Let $M_i = \Gamma_i \backslash \mathbb{H}$ $(i \in I)$ be a family of length-commensurable Riemann surfaces, with $\Gamma_i \subset \mathrm{SL}_2(\mathbb{R})$ Zariski-dense. Then the associated quaternion algebras $A_{\Gamma_i}$ $(i \in I)$ belong to \emph{finitely many} isomorphism classes (over the common center).
\end{thm}

This theorem is only one application of the investigation into the {\it genus problem} for algebraic groups and the finiteness conjecture for  the algebraic hulls of weakly commensurable Zariski-dense subgroups that relies on the consideration of good reduction (cf. \cite{R-ICM}, \cite[\S\S 8, 9]{RR-Survey}). More geometric applications to locally symmetric spaces (including those of infinite volume) can potentially be given using the approach developed in \cite{PrR}. This provides  rather encouraging examples of the use in differential geometry of new techniques based on good reduction that are likely to be helpful in future research.

\vskip2mm

\noindent {\small {\bf Acknowledgements.} We would like to thank Skip Garibaldi, Danny Krashen, Steven Sam, Sergey Tikhonov, Charlotte Ure, and the anonymous referees for their comments and suggestions that helped to improve the exposition.}


\bibliographystyle{amsplain}

\end{document}